\documentclass[11pt]{amsart}

\usepackage{amscd}
\usepackage{amsthm, amsfonts, amssymb}
\usepackage
{hyperref}

\theoremstyle{definition}
\newtheorem{ntn}{Notation}[section]

\theoremstyle{plain}
\newtheorem{lem}[ntn]{Lemma}
\newtheorem{prp}[ntn]{Proposition}
\newtheorem{thm}[ntn]{Theorem}
\newtheorem{cor}[ntn]{Corollary}

\theoremstyle{remark}

\newtheorem{rem}[ntn]{Remark}
\newtheorem{exa}[ntn]{Example}

\numberwithin{equation}{section}

\newcommand{\N}{\mathbb{N}}
\newcommand{\z}{\mathbb{Z}}
\newcommand{\q}{\mathbb{Q}}
\newcommand{\R}{\mathbb{R}}
\newcommand{\C}{\mathbb{C}}
\newcommand{\F}{\mathbb{F}}

\newcommand{\ppp}{\mathfrak{p}}
\newcommand{\mmm}{\mathfrak{m}}

\newcommand{\GL}{\mathit{{\rm GL}}}
\newcommand{\SL}{\mathit{{\rm SL}}}
\newcommand{\GM}{\mathit{{\rm GM}}}
\newcommand{\SK}{\mathit{{\rm SK}}}
\newcommand{\Ee}{\mathit{{ \rm E}}}

\renewcommand{\H}{\tilde{H}}
\newcommand{\rr}{{R^\times}}
\newcommand{\fff}{{R^\times}}

\newcommand{\tors}{{{\rm Tor}_1^{\z}}}

\newcommand{\exts }{{{\rm Ext}_{\z}^1}}

\newcommand{\mt}{\mapsto}
\newcommand{\lan}{\langle}
\newcommand{\ran}{\rangle}

\newcommand{\se}{\subseteq}
\newcommand{\arr}{\rightarrow}
\newcommand{\larr}{\longrightarrow}

\newcommand{\two}{\twoheadrightarrow}

\newcommand{\id}{{\rm id}}
\newcommand{\im}{{\rm im}}
\newcommand{\ind}{{\rm ind}}

\newcommand{\inc}{{\rm inc}}

\newcommand{\diag}{{\rm diag}}
\newcommand{\coker}{{\rm coker}}
\renewcommand{\char}{{\rm char}}

\newcommand {\mtx}[4]
{\left(
\begin{array}{cc}
#1 & #2   \\
#3 & #4
\end{array}
\right)}

\newtheoremstyle{athm}
  {}
  {}
  {\itshape}
  {}
  {\scshape}
  {}
  {.5em}
  {\thmnote{#3}}
\theoremstyle{athm}
\newtheorem*{athm}{}

\begin{document}
\title[Bloch-Wigner theorem]{Bloch-Wigner theorem over rings with many units II}
\author{Behrooz Mirzaii}
\author{Fatemeh Y. Mokari}

\begin{abstract}
In this article we prove a generalization of the Bloch-Wigner exact sequence over 
commutative rings with many units. When the ring is a domain, we get a generalization 
of Suslin's Bloch-Wigner exact sequence over infinite fields. Our proof is different and is 
easier than Suslin's proof, even in its general form. But nevertheless we use some of 
Suslin's results which relates the Bloch group of the ring to the third homology group 
of the general linear group of the ring. From there we take an easier path.  
\end{abstract}
\maketitle

%%%%%%%%%%%%%%%%%%%%%%%%%%%%%%%%%%%%%%%%%%%%%%%%%%%%%%%%%%%%%%%%%%%%%%%%%%%%%%%%
\section*{Introduction}
%%%%%%%%%%%%%%%%%%%%%%%%%%%%%%%%%%%%%%%%%%%%%%%%%%%%%%%%%%%%%%%%%%%%%%%%%%%%%%%%

For a commutative ring $R$ with 1, there are two types of $K$-groups: Milnor
$K$-groups and Quillen $K$-groups, denoted by $K_n^M(R)$ and  $K_n(R)$, respectively.
For any positive integer $n$, there is a canonical homomorphism 
\[
\iota_n: K_n^M(R)\arr K_n(R).
\]
When $R$ is an infinite field or more generally a ring with many units,
$\iota_1$ and $\iota_2$ are isomorphisms  \cite{vdkallen1977}.
This fails for $\iota_n$ when $n \geq 3$ most of the times. 

Let $R$ be a commutative ring with many units. Then it is well-known that the kernel of $\iota_n$ 
is annihilated by multiplication  by $(n-1)!$ \cite{suslin1985}, \cite{nes-suslin1990},
\cite{guin1989}, and the cokernel of $\iota_n$ usually is very large. The group
\[
K_n(R)^\ind:=\coker(K_n^M(R)\overset{\iota_n}{\larr} K_n(R))
\]
is called the indecomposable part of $K_n(R)$. Clearly $K_n(R)^\ind=0$ for $n=1,2$. 
Although there are nice structural theorems for $K_n(R)^\ind$,
but still there are many unanswered conjectures related to its structure \cite{soule1985}. 

The original Bloch-Wigner exact sequence, proved by Bloch and Wigner in somewhat different 
form but not published by them, asserts the existence of the exact sequence
\[
\begin{array}{c}
0 \arr \q/\z \arr H_3(\SL_2(\C),\z) \arr \ppp(\C) \arr \bigwedge_\z^2 \C^\times \arr K_2(\C) \arr 0,
\end{array}
\]
where $\ppp(\C)$ is the pre-Bloch group of $\C$ and the map $\ppp(\C) \arr \bigwedge_\z^2 \C^\times$
is given by $[a]\mt a \wedge (1-a)$ \cite{bloch2000}.  One can show that
$H_3(\SL_2(\C),\z)\simeq K_3^\ind(\C)$
and so the Bloch-Wigner exact sequence gives a precise description of the 
indecomposable part of $K_3(\C)$. These facts easily can be generalized to all algebraically 
closed field of characteristic zero \cite{dupont-sah1982}, \cite{sah1989}. 

The Bloch-Wigner exact sequence appears in different areas of mathematics, such as 
number theory \cite{bloch2000}, three dimensional hyperbolic geometry 
\cite{dupont-sah1982}, algebraic $K$-theory \cite{suslin1991}, etc. 
Therefore it is very desirable to have a nice generalization of this exact sequence 
to other class of rings.

In a remarkable paper \cite{suslin1991}, Suslin has generalized this result to all infinite
fields. In fact he proved that for any infinite field $F$, we have the exact sequence
\begin{align*}
0 \arr \!\tors(\mu(F), \mu(F))^\sim \!\arr \!K_3(F)^\ind \!\arr
\ppp(F)\! \arr (F^\times\otimes_\z F^\times)_\sigma\!\arr K_2(F) \arr 0,
\end{align*}
where $\tors(\mu(F), \mu(F))^\sim$ is the unique non-trivial extension of the group 
$\tors(\mu(F), \mu(F))$ by $\z/2$. 

In present article we generalize the Bloch-Wigner exact sequence over rings 
with many units. Here is our main theorem.

\begin{athm}[{\bf Theorem \ref{mirzaii-mokari}.}]
Let $R$ be a commutative ring with many units. Then we have the exact sequence
\begin{align*}
T_R  \arr K_3^\ind(R) \arr \ppp(R) \arr  (\rr\otimes_\z \rr)_\sigma \arr K_2(R) 
\arr 0,
\end{align*}
where $T_R$ sits in an exact sequence
\[
0 \arr \tors(\mu(R),\mu(R))_\sigma \arr T_R \arr
H_1(\Sigma_2, \mu_{2^\infty}(R)\otimes_\z \mu_{2^\infty}(R)) \arr 0.
\]
Moreover when $R$ is an integral domain, we have the exact sequence
\begin{align*}
0 \arr \tors(\mu(R), \mu(R))^\sim\! \arr K_3^\ind(R)\! \arr \ppp(R)
\arr \! (\rr\otimes_\z \rr)_\sigma\! \arr K_2(R) \arr 0,
\end{align*}
where the composition $\tors(\mu(R),\mu(R)) \arr \tors(\mu(R), \mu(R))^\sim \arr K_3^\ind(R)$ is 
induced by the map $\mu(R) \arr \SL_2(R)$,  $\xi \mt \diag(\xi, \xi^{-1})$.
\end{athm}

Our proof of this theorem is new and easier than Suslin's proof. But nevertheless
we use a result of Suslin which relates the Bloch group of $R$ to the third 
homology group of $\GL_3(R)$ (see Theorem \ref{gl-gm3} below). 

There are two new, but simple, ingredients in our proof. The first one is the explicit 
description of the composition map
\[
K_n^M(R) \overset{\iota_n}{\larr} K_n(R) \overset{h_n}{\larr} H_n(\GL(R),\z),
\]
where $h_n$ is the Hurewicz map (see Section \ref{section1}). In fact this follows 
easily from another result of Suslin (see Proposition \ref{hur1} below). But it seems that, 
except for lower $K$-groups, such an explicit formula did not appear anywhere in the 
literature. (But nevertheless see \cite[Section~4]{mirzaii2008}.) When $n=3$, this allows 
us to give a simple description of $K_3^\ind(R)$ as a quotient of $H_3(\SL(R), \z)$. 

The second idea comes from the observation that for any integral domain $R$ of $\char(R)\neq 2$
and with finite $\mu_{2^\infty}(R)$, $T_R$ is the unique nontrivial extension of $\z/2$ by  
$\tors(\mu(R),\mu(R))$ and we show that this is isomorphic to $\tors(\mu(R),\mu(R))^\sim$.

As an application we show that if $R$ is a discrete valuation ring with field of fraction $K$
and infinite residue field $F$ such that $\char(F)=\char(K)$, then 
\[
B(R)\simeq B(K),
\]
where $B(R)$ is the Bloch group of $R$. This improves a result of Hutchinson which 
proves the same isomorphism, but up to 2-torsion elements \cite{hutchinson2014}.
Furthermore we will study the Bloch-Wigner exact sequence for rings of dual
numbers $R[\epsilon]$, where $R$ is a domain with many units, which is of some 
interest \cite{cathelineau2007}. 
At the end we use these information to study the relation between $K_3^\ind(R)$ and
the third homology group of $\SL_2(R)$.

In this paper by $H_n(G)$ we mean  the $n$-th homology of the group $G$ with integral coefficients, 
namely $H_n(G, \z)$. If $A \arr A'$ is a homomorphism of abelian groups, by $A'/A$ we
mean $\coker(A \arr A')$ and we take other liberties of this kind. For a group $A$, by 
$A_{\z[1/2]}$ we mean $A\otimes_\z \z[1/2]$. Here all rings are commutative and have $1$.

%%%%%%%%%%%%%%%%%%%%%%%%%%%%%%%%%%%%%%%%%%%%%%%%%%%%%%%%%%%%%%%%%%%%%%%%%%%
\section{{\it K}-theory and homology of general linear groups}\label{section1}
%%%%%%%%%%%%%%%%%%%%%%%%%%%%%%%%%%%%%%%%%%%%%%%%%%%%%%%%%%%%%%%%%%%%%%%%%%%

For a commutative ring $R$ and an integer $n$, let $K_n(R):=\pi_n(B\GL(R)^+)$
be the $n$-th Quillen $K$-group, where $B\GL(R)^+$ is the plus-construction
of the classifying space of $\GL(R)$ with respect to the elementary subgroup $\Ee(R)$
\cite[Chap.~I]{loday1976}. If $n\ge 2$, then 
$K_n(R)\simeq \pi_n(B\Ee(R)^+)$, where $B\Ee(R)^+$ is the plus-construction of $B\Ee(R)$ 
with respect to $\Ee(R)$ \cite[Proposition~1.1.7]{loday1976}. On the other hand,
$H_n(B\Ee(R)^+,\z)\simeq H_n(\Ee(R))$ and $H_n(B\GL(R)^+,\z)\simeq H_n(\GL(R))$,
thus for $n \geq 2$, the Hurewicz homomorphisms induce the commutative diagram
\[
\begin{array}{ccc}
K_n(R) \!\!\!& \overset{h_n}{-\!\!\!-\!\!\!-\!\!\!-\!\!\!-\!\!\!-\!\!\!
-\!\!\!\larr}& \!\!\!\!\! H_n(\GL(R))\\
{}_{h_n'}\!\!\!\searrow\!\!\!\!\!\!\!\!\!\!\!\!\!\!\!\!\! & &\!\!\!
\hspace{-2.2cm} \nearrow\\ 
& H_n(\Ee(R)). &
\end{array}
\]
Moreover the space $B\GL(R)^+$ has a natural $H$-space structure. This allows us to construct 
product maps 
\[
K_m(R) \otimes_\z K_n(R) \arr K_{m+n}(R)
\] 
which is anti-commutative \cite[Chap.~II]{loday1976}.
It is easy to see that $K_1(R)\overset{h_1}{\simeq} H_1(\GL(R))$ and it is well-known 
that $K_2(R)\overset{h_2'}{\simeq} H_2(\Ee(R))$. Moreover 
$h_3':K_3(R) \arr H_3(\Ee(R))$ is surjective and its kernel coincides with the image of 
the composition map
$K_1(\z)\otimes_\z K_2(R) \arr K_1(R) \otimes_\z K_2(R)\arr K_3(R)$ 
\cite[Corollary~5.2]{suslin1991}, which we denote this image by $l(-1)K_2(R)$.
Note that $K_1(\z)\simeq \{-1,1\}$.
Thus we have the exact sequence
\begin{equation}\label{k3-exact}
K_1(\z)\otimes_\z K_2(R) \arr K_3(R) \arr H_3(\Ee(R)) \arr 0.
\end{equation}

The $n$-th Milnor $K$-group of a commutative ring $R$
is defined as abelian group $K_n^M(R)$ generated by symbols
$\{a_1, \dots, a_n\}$, $a_i \in \rr$,  $i=1, \dots, n$,
subject to the following relations
\medskip
\par (i) $\{a_1, \dots,a_ia_i', \dots, a_n\!\}\!= \!\{a_1, \dots,a_i, \dots, a_n\!\}\!+\!
\{a_1, \dots,a_i', \dots, a_n\!\}$, any $i$,
\medskip
\par (ii) $\{a_1,  \dots, a_n\}=0$ if there exist $i,j$, $ i \neq j$,
such  that $a_i+a_j=0$ or $1$.
\medskip
~\\
Clearly we have the anti-commutative product map
\[
K_m^M(R) \otimes_\z K_n^M(R) \arr K_{m+n}^M(R),
\]
\[ 
\{a_1, \dots,a_m\}\otimes \{b_1, \dots, b_n\} \mt \{a_1, \dots,a_m, b_1, \dots, b_m\}.
\]
Using the product structures for Quillen and Milnor $K$-groups we obtain a 
natural map 
\[
\iota_n: K_n^M(R) \arr K_n(R).
\]

We say that a commutative ring $R$ is a ring with many units if 
for any finite number of surjective linear
forms $f_i: R^2 \arr R$, there exists a $v \in  R^2$ such that, for
all $i$, $f_i(v) \in \rr$. Important examples of rings with many units
are semilocal rings with infinite residue fields. For more on these rings
we refer the reader to \cite{vdkallen1977}, \cite{guin1989}
or \cite{mirzaii2008}.

Over these rings, $K_1^M(R) =\rr \simeq K_1(R)$.
Since $K_1(R) \simeq \rr \times \SK_1(R)$, we have
$\SK_1(R):=\SL(R)/\Ee(R)=0$ and hence $\Ee(R)=\SL(R)$.
Moreover by a theorem of  Van der Kallen \cite[Theorem~8.4]{vdkallen1977} we have 
the isomorphisms
\[
 K_2^M(R)\simeq \rr \otimes_\z \rr/\lan a \otimes
(1-a):a, 1-a \in \rr \ran \simeq K_2(R),
\]
which is given by $\{a,b\}\mt {\bf c}(\diag(a,1,a^{-1}), \diag(b,b^{-1},1))$.

For an arbitrary group $G$, let $B_\bullet(G)\overset{\varepsilon}{\arr}\z$
denote the bar resolution of $G$. We turn the $B_n(G)$ into a right $G$-module in usual way.
For any left $G$-module $N$, the homology group $H_n(G, N)$ coincides with the homology
of the complex $B_\bullet(G) \otimes_G N$. In particular
\[
H_n(G)=H_n(B_\bullet(G) \otimes_G \z)=H_n(B_\bullet(G)_G).
\]
For simplicity the element of $B_n(G)$ represented by $[g_1|\dots|g_n]$ again is denoted
by $[g_1|\dots|g_n]$. Let
\[
{\rm \bf{c}}({g}_1, {g}_2,\dots, {g}_n):=\sum_{\sigma \in \Sigma_n} {{\rm
sign}(\sigma)}[{g}_{\sigma(1)}| {g}_{\sigma(2)}|\dots|{g}_{\sigma(n)}] \in
H_n(G),
\]
where ${g}_i \in G$ pairwise commute and $\Sigma_n$ is the symmetric
group of degree $n$. The proof of the following facts are straightforward: 
\par {\rm (i)} If $h_1\in G$ commutes with all the elements
$g_1, \dots, g_n \in G$, then
\[
{\rm \bf{c}}(g_1h_1, g_2,\dots, g_n)= {\rm \bf{c}}(g_1, g_2,\dots,
g_n)+{\rm \bf{c}}(h_1, g_2,\dots, g_n),
\]
\par {\rm (ii)}
For every $\sigma \in \Sigma_n$, ${\rm \bf{c}}(g_{\sigma(1)},\dots,
g_{\sigma(n)})={\rm sign(\sigma)} {\rm \bf{c}}(g_1,\dots,g_n)$,
\par {\rm (iii)}
The cup product of ${\rm \bf{c}}(g_1,\dots, g_p)\!\in\!\! H_p(G)$
and ${\rm \bf{c}}(g_1',\dots, g_q')\! \in \!\!H_q(G')$ is 
$${\rm
\bf{c}}((g_1, 1), \dots, (g_p,1),(1,g_1'), \dots, (1,g_q')) \in
H_{p+q}(G \times G').$$

If $R$ is a commutative ring with many units, then we have the homological stability
isomorphisms
\[
H_n(\GL_n(R)) \overset{\sim}{\larr}  H_n(\GL_{n+1}(R)) \overset{\sim}{\larr}
H_n(\GL_{n+2}(R)) \overset{\sim}{\larr} \cdots.
\]
Moreover one can construct a surjective map $H_n(\GL_{n}(R)) \two  K_n^M(R)$ 
such that the sequence
\[
H_n(\GL_{n-1}(R)) \arr H_n(\GL_{n}(R)) \arr  K_n^M(R) \arr 0
\]
is exact and the isomorphism 
\[
K_n^M(R) \overset{\simeq}{\larr} H_n(\GL_{n}(R))/H_n(\GL_{n-1}(R))
\]
is given by the formula $\{a_1, \dots, a_n\} \mt a_1\cup \dots \cup a_n 
\mod \im H_n(\GL_{n-1}(R))$
\cite{nes-suslin1990}, \cite{suslin1985}, \cite{guin1989}. 

%%%%%%%%%%%%%%%%%%%%%%%%%%%%%%%%%%%%%%%%%%%%%%%%%%%%%%%%%%%%%%%%%%%%%%%%%%%
\section{Milnor {\it K}-groups and the hurewicz homomorphism}
%%%%%%%%%%%%%%%%%%%%%%%%%%%%%%%%%%%%%%%%%%%%%%%%%%%%%%%%%%%%%%%%%%%%%%%%%%%

Let $R$ be a ring with many units.
In the previous section we gave an explicit description of the map 
\[
l_n:=h_n\circ \iota_n: K_n^M(R) \arr H_n(\GL(R))\simeq H_n(\GL_n(R))
\]
for $n=2$. In fact we have
\[
l_2(\{a,b\})={\bf c}(\diag(a,1,a^{-1}), \diag(b,b^{-1},1))
={\bf c}(\diag(a,1), \diag(b,b^{-1})).
\]
We will generalize this formula to any $n$. 
In fact this follows from the following result of Suslin. 

\begin{prp}\label{hur1}
Let $R$ be a ring with many units and let $x \in K_m(R)$ and $y \in K_n(R)$. Then
\[
h_{m+n}(x.y)={\psi_{m,n}}_\ast(h_m(x) \times h_n(y))-mn (h_m(x)\cup h_n(y)),
\]
where $\psi_{m,n}:\GL_m(R) \times \GL_n(R)\arr \GL_{mn}(R)$ is the matrix tensor product and
the cup product is induced by 
\[
\GL_m(R) \times \GL_n(R)\arr \GL_{m+n}(R),\ \ \ (A,B) \mt {\mtx A {0} {0} B}.
\]
\end{prp}
\begin{proof}
Suslin has proved this for infinite fields \cite[Lemma 4.2]{suslin1985}. The same 
proof works, without any changes, for rings with many units. It should be mentioned 
that for the proof we need the homological stability  isomorphisms. Here
${\psi_{m,n}}_\ast(h_m(x) \times h_n(y))$ is the image of
$h_m(x) \otimes h_n(y)$ under the composition  
\begin{align*}
H_m(\GL_m(R)) \otimes_\z H_n(\GL_n(R)) &\arr H_{m+n}(\GL_{m}(R) \times \GL_{n}(R))\\
& \arr H_{m+n}(\GL_{mn}(R))\\
& \arr H_{m+n}(\GL(R))\simeq H_{m+n}(\GL_{m+n}(R)).
\end{align*}
\end{proof}

The following results is very important for us.

\begin{cor}\label{hur2}
Let $R$ be a ring with many units. Then the homomorphism
$l_n: K_n^M(R) \arr H_n(\GL_n(R))$ is given by the formula
\[
l_n(\{a_1, \dots, a_n\})=[a_1, \dots, a_n]:= {\bf c}(A_{1, n}, \dots , A_{n, n}),
\]
where $A_{i, n}:=\diag(a_i, \dots, a_i, a_i^{-(i-1)},I_{n-i}) \in \GL_n(R)$.
\end{cor}
\begin{proof}
We prove this by induction on $n$. The case $n=1$ is clear. By induction assume
that $l_n(\{a_1, \dots, a_{n-1}\})=[a_1, \dots, a_{n-1}]$. Since the map
$\psi_{n-1,1}: \GL_{n-1}(R) \times \GL_{1}(R) \arr \GL_{n-1}(R)$ is given by $(A,a) \mt aA$,
by Proposition~\ref{hur1}, we have
\begin{align*}
l_n(\{a_1, \dots, a_n\})&=h_n(\{a_1, \dots, a_{n-1}\}.\{a_n\})\\
& = {\psi_{n-1,1}}_\ast\bigg({\bf c}(A_{1, n-1}, \dots , A_{n-1, n-1}) \times {\bf c}(a_n)\bigg)\\
& - (n-1) {\bf c}(A_{1, n-1}, \dots , A_{n-1, n-1}) \cup {\bf c}(a_n)\\
& = {\bf c}(A_{1, n-1}, \dots , A_{n-1, n-1}, a_nI_{n-1})\\
& - (n-1){\bf c}(A_{1, n}, \dots , A_{n-1, n}, \diag(I_{n-1},a_n))\\
& ={\bf c}(A_{1, n}, \dots , A_{n-1, n}, \diag(a_nI_{n-1},a_n^{-(n-1)}))\\
& ={\bf c}(A_{1, n}, \dots , A_{n, n})\\
& = [a_1, \dots, a_n].
\end{align*}
\end{proof}

\begin{rem}
The above proposition shows that the map $l_n$ coincides with the map
$\nu_n:K_n^M(R) \arr H_n(\GL_n(R))$ in \cite[Section~4]{mirzaii2008}, 
which is constructed directly.
\end{rem}

It is very easy to show that $[a_1, \dots, a_n] \in \!H_n(\SL(R))$. For example
in $H_3(\GL_3(R))$ we have
${\rm \bf{c}}(\diag(a,1,a^{-1}),\diag(b, b^{-1},1),\diag(1, c^{-1}, c))=0$,
and thus
\[
[a, b,c]={\rm \bf{c}}(\diag(a,1,a^{-1}),
\diag(b, b^{-1},1),\diag(c, 1, c^{-1})) \in H_3(\SL(R)).
\]
\begin{cor}\label{k3-h3}
If $R$ is a ring with many units, then the Hurewicz map
$h_3':K_3(R) \arr H_3(\SL(R))$
induces an isomorphism $K_3^\ind(R) \simeq H_3(\SL(R))/T$, 
where $T=\lan [a,b,c]:a,b,c \in \rr \ran$.
\end{cor}
\begin{proof}
By Corollary \ref{hur2}, $\overline{h}'_3:K_3^\ind(R) \arr H_3(\SL(R))/T$
is well-defined and by (\ref{k3-exact}) it is surjective. 
The injectivity of this map follows from the fact that
$\ker(K_3(R) \two H_3(\SL(R)))\se \im(K_3^M(R) \arr K_3(R))$,
and this follows from exact sequence (\ref{k3-exact}).
\end{proof}

%%%%%%%%%%%%%%%%%%%%%%%%%%%%%%%%%%%%%%%%%%%%%%%%%%%%%%%%%%%%%%%%%%%%%%%%%%%%%
\section{Bloch group and homology of general liner groups}
%%%%%%%%%%%%%%%%%%%%%%%%%%%%%%%%%%%%%%%%%%%%%%%%%%%%%%%%%%%%%%%%%%%%%%%%%%%%%

Define the pre-Bloch group $\ppp(R)$ of a commutative ring $R$ as the quotient of
the free abelian group $Q(R)$ generated by symbols $[a]$, $a, 1-a \in \fff$,
by the subgroup generated by elements of the form
\[
[a] -[b]+\bigg[\frac{b}{a}\bigg]-\bigg[\frac{1- a^{-1}}{1- b^{-1}}\bigg]
+ \bigg[\frac{1-a}{1-b}\bigg],
\]
where $a, 1-a, b, 1-b, a-b  \in \fff$. Let
\[
\lambda': Q(R) \arr \fff \otimes_\z \fff, \ \ \ \ [a] \mapsto a \otimes (1-a).
\]
By a direct computation, we have
\[
\lambda'\Big(
[a] -[b]+\bigg[\frac{b}{a}\bigg]-\bigg[\frac{1- a^{-1}}{1- b^{-1}}\bigg]
+ \bigg[\frac{1-a}{1-b}\bigg] \Big)
=a \otimes \bigg( \frac{1-a}{1-b}\bigg)+\bigg(\frac{1-a}{1-b}\bigg)\otimes a.
\]
Let 
\[
(\fff \otimes_\z \fff)_\sigma :=\fff \otimes_\z \fff/
\lan a\otimes b + b\otimes a: a, b \in \fff \rangle
\]
(see the next section for the choice of the notation).
We denote the elements of $\ppp(R)$ and $(\fff \otimes_\z \fff)_\sigma$
represented by $[a]$ and $a \otimes b$ again by $[a]$ and $a\otimes b$,
respectively. Thus we have a well-defined homomorphism
\[
\lambda: \ppp(R) \arr (\fff \otimes_\z \fff)_\sigma, \ \ \
[a] \mapsto a \otimes (1-a).
\]
The kernel of $\lambda$ is called the Bloch group of $R$ and is denoted by $B(R)$.
If $R$ is a ring with many units, then we obtain the exact sequence
\[
0 \arr B(R) \arr \ppp(R) \arr (\fff \otimes_\z \fff)_\sigma \arr K_2^M(R) \arr 0.
\]
The following result is due to Suslin \cite[Theorem 2.1, Proposition~3.1]{suslin1991}.

\begin{thm}\label{gl-gm3}
Let $R$ be a ring with many units. Then we have the exact sequences
\[
H_3(\GM_2(R)) \arr H_3(\GL_2(R)) \arr B(R) \arr 0,
\]
\[
H_3(\GM_2(R))\oplus H_3(T_3) \arr H_3(\GL_3(R)) \arr B(R) \arr 0,
\]
where $\GM_2(R)$ is the subgroup of monomial matrices in $\GL_2(R)$ and 
$T_3$ is the subgroup of diagonal matrices in $\GL_3(R)$.
Moreover the homomorphism $H_3(\GL_2(R)) \arr B(R)$ factors through
$H_3(\GL_3(R)) \arr B(R)$.
\end{thm}
\begin{proof}
Suslin has proved these results over infinite fields 
\cite[Theorem~2.1, Proposition~3.1]{suslin1991}.
But his arguments also work in our situation. 
\end{proof}

%%%%%%%%%%%%%%%%%%%%%%%%%%%%%%%%%%%%%%%%%%%%%%%%%%%%%%%%%%%%%%%%%%%%%%%%%%%
\section{Third homology of $\GM_2$}
%%%%%%%%%%%%%%%%%%%%%%%%%%%%%%%%%%%%%%%%%%%%%%%%%%%%%%%%%%%%%%%%%%%%%%%%%%%

Let $T_2:=\rr \times \rr$ be the diagonal subgroup of $\GM_2(R)$. Then 
$\GM_2(R) = T_2 \rtimes \Sigma_2$,
where $\Sigma_2=\Bigg\{ {\mtx 1 0 0 1}, {\mtx 0 1 1 0} \Bigg\}$. We often think of
$\Sigma_2$ as the symmetric group of order two $\{1, \sigma\}$, where the action 
of $\sigma$ on $T_2$ is given by $\sigma(a,b)=(b,a)$.
In this section we will study the associated Lyndon-Hochschild-Serre spectral 
sequence of the extension $1 \arr T_2 \arr \GM_2(R) \arr \Sigma_2 \arr 1$;
\[
E^2_{p,q}=H_p(\Sigma_2 ,H_q(T_2)) \Rightarrow H_{p+q}(\GM_2(R)).
\]
For a $\Sigma_2$-module $M$, we have 
\[
H_p(\Sigma_2, M)\simeq
\begin{cases}
M_{\Sigma_2} & \text{if $p=0$} \\
M^{\Sigma_2}/(1+\sigma)(M) & \text{if $p$ is odd} \\
\ker(M_{\Sigma_2}\overset{1 +\sigma}{-\!\!\!-\!\!\!\larr} M^{\Sigma_2})   &   
\text{if $p$ is even.}
\end{cases}
\]
For simplicity we denote $M_{\Sigma_2}$ and $M^{\Sigma_2}$ by $M_{\sigma}$ and $M^{\sigma}$,
respectively.

\begin{lem}\label{h2}
Let $A$ be an abelian group and let $\Sigma_2$ acts on $A\oplus A$ and $A\otimes_\z A$ as
$\sigma(a,b)=(b,a)$ and $\sigma(a\otimes b)=-b\otimes a$, respectively. Then 
\par {\rm (i)} $H_0(\Sigma_2,  A\oplus A)\simeq A$ and 
$H_0(\Sigma_2,  A\otimes_\z A) \simeq (A\otimes_\z A)_\sigma$.
\par {\rm (ii)} For any $p \ge 1$, $H_p(\Sigma_2,  A\oplus A)=0$. 
\par {\rm (iii)} If $A\simeq \bigoplus_{i=1}^n A_i$, then 
$H_p(\Sigma_2,  A\otimes_\z A)\simeq \bigoplus_{i=1}^n H_p(\Sigma_2,  A_i\otimes_\z A_i)$.
\par {\rm (iv)} If $p$ is odd, then 
$H_p(\Sigma_2,  A\otimes_\z A)\simeq H_p(\Sigma_2, {}_{2^\infty}\!A\otimes_\z {}_{2^\infty}\!A)$.
\par {\rm (v)} If $p$ is even, then 
$H_p(\Sigma_2,  A\otimes_\z A)=
\lan \overline{a \otimes a}\in (A\otimes_\z A)_\sigma: a\in A\ran$.
Here ${}_{2^\infty}\!A$ is the subgroup of $2$-torsion elements of $A$, i.e.
\[
{}_{2^\infty}\!A:=\Big\{a \in A: there \ exists\ an \   m \in \N,\  s.t.\  {2^m}a=0\Big\}.
\]
\end{lem}
\begin{proof}
Parts (i) and (ii) follow from a direct and easy computation. In part (iii) 
we may assume that $n=2$. Then 
\begin{align*}
H_p(\Sigma_2,A\otimes_\z A)& \simeq H_p(\Sigma_2, A_1\otimes_\z A_1) 
\oplus H_p(\Sigma_2, A_2\otimes_\z A_2)\\
& \oplus H_p(\Sigma_2, A_1\otimes_\z A_2 \oplus A_2\otimes_\z A_1),
\end{align*}
where the action of $\sigma$ on $A_1\otimes_\z A_2 \oplus A_2\otimes_\z A_1$ is given by
\[
\sigma(a_1\otimes a_2, b_2\otimes b_1)=-(b_1 \otimes b_2, a_2\otimes a_1).
\] 
Now as part (ii), $H_p(\Sigma_2, A_1\otimes_\z A_2 \oplus A_2\otimes_\z A_1)=0$.
Let $\{A_j:j \in J\}$ be a family of finitely generated subgroups of $A$ such that
$J$ is a directed set and $A\simeq \underset{\underset{J}{\larr}}{\lim} A_j$.
Then $A \otimes_\z A \simeq \underset{\underset{J}{\larr}}{\lim}( A_j \otimes_\z A_j)$. 
Thus in parts 
(iv) and (v) we may assume that $A$ is finitely generated and by part (iii)
we may assume that $A$ is cyclic, which in this case the claim is easy to prove.
\end{proof}

By Lemma \ref{h2}, some of the $E^2$-terms of the spectral sequence can be computed as follow:
\[
E^2_{p,0}\!\simeq\!
\begin{cases} \z & \text{if $p=0$} \\
\z/2 &   \text{if $p$ is odd}      \\
0    &   \text{if $p$ is even}
\end{cases},
E^2_{p,1}
\!\simeq\!
\begin{cases}
\fff & \text{if $p=0$} \\
0    &   \text{if $p \neq 0$}
\end{cases},
E_{p,2}^2 \!\simeq \!H_p(\Sigma_2, \rr \otimes_\z \rr).
\]
From these we obtain the isomorphism $E_{1,2}^\infty \simeq E_{1,2}^2$. Moreover
by Lemma~\ref{h2}, $E_{2,2}^2$ is generated by the elements $\overline{a \otimes a}$, $a\in \rr$,
and a direct computation shows that $d_{2,2}^2(\overline{a \otimes a})=0$.
This implies that 
$E_{0,3}^\infty\simeq E_{0,3}^2 = H_3(T_2)_\sigma$.
The spectral sequence $E_{p,q}^1$ gives us a filtration
\[
0=F_{-1}H_3(\GM_2(R)) \se \cdots \se F_3H_3(\GM_2(R))=H_3(\GM_2(R)),
\]
such that
\[
\begin{array}{l}
\vspace{1.5 mm}
E^\infty_{0, 3}\simeq F_0H_3(\GM_2(R))=H_3(T_2)_\sigma,\\
\vspace{1.5 mm}
E^\infty_{1, 2} \simeq F_1H_3(\GM_2(R))/F_0H_3(\GM_2(R))\simeq E_{1,2}^2,\\
\vspace{1.5 mm}
E^\infty_{2, 1}\simeq F_2H_3(\GM_2(R))/F_1H_3(\GM_2(R))=0,\\
E^\infty_{3, 0}\simeq 
H_3(\GM_2(R))/F_2H_3(\GM_2(R))\simeq H_3(\Sigma_2).
\end{array}
\]
Thus
\begin{align}\label{GM2}
H_3(\GM_2(R)) \simeq F_2H_3(\GM_2(R)) \oplus H_3(\Sigma_2),
\end{align}
\begin{align}
E_{1, 2}^2 \simeq F_2H_3(\GM_2(R))/H_3(T_2)_\sigma.
\end{align}
Set $M:=H_3(\rr) \oplus H_3(\rr) \oplus \rr \otimes_\z H_2(\rr)
\oplus H_2(\rr) \otimes_\z \rr \se H_3(T_2)$.
By applying the Snake lemma to the commutative diagram
\[
\begin{CD}
& M_\sigma &\hspace{-0.7 cm}{-\!\!\!-\!\!\!-\!\!\!-\!\!\!-\!\!\!-\!\!\!-\!\!\!\arr}\!\!\!\!\!\! &
H_3(T_2)_\sigma & \!\!\!\!\!\!-\!\!\!-\!\!\!\arr  & \tors(\mu(R),\mu(R))_\sigma & \arr 0 \\
 &  @VVV  @VVV  @VVV & \\
0 \arr & F_2H_3(\GM_2(R)) & \arr & F_2H_3(\GM_2(R)) &
{-\!\!\!-\!\!\!-\!\!\!-\!\!\!-\!\!\!-\!\!\!\arr}\hspace{-1.4 cm} & 0 &
\hspace{-1.3 cm} -\!\!\!-\!\!\!-\!\!\!-\!\!\!-\!\!\!-\!\!\!-\!\!\!\arr 0
\end{CD}
\]
and Lemma \ref{h2} we obtain the exact sequence
\begin{equation}\label{TR-exact}
0 \arr \tors(\mu(R),\mu(R))_\sigma \arr T_R \arr
H_1(\Sigma_2, \mu_{2^\infty}(R) \otimes_\z \mu_{2^\infty}(R)) \arr 0,
\end{equation}
where $T_R:=F_2H_3(\GM_2(R))/M_\sigma=F_2H_3(\GM_2(R))/M$.
To give an explicit description of the composition map
\[
\begin{CD}
H_2(T_2)^\sigma -\!\!\!\twoheadrightarrow E_{1,2}^2 \overset{\simeq}{\larr}
F_2H_3(\GM_2(R))/H_3(T_2) \se H_3(\GM_2(R))/H_3(T_2),
\end{CD}
\]
we need to introduce certain notations.

Let $G$ be a group and let $M$ be a $G$-module. Any element $g \in G$, determines an 
automorphism of the complex $B_\bullet(G) \otimes_G N$ given by
\[
[g_1|\dots|g_n]\otimes m \mapsto [gg _1 g^{-1}|\dots|gg_ng^{-1}]\otimes gm.
\]
This automorphism is homotopic to the identity, with the
corresponding homotopy given by the formula
\[
\rho_g: B_n(G) \otimes_G N   \arr     B_{n+1}(G) \otimes_G N ,
\]
\[
[g_1|\dots|g_n]\otimes m \mapsto  \sum_{j=0}^n (-1)^j
[g_1|\dots|g_j|g^{-1}|g g_{j+1} g^{-1}|\dots|g g_n g^{-1}] \otimes m.
\]
The following lemma is similar to \cite[Lemma 2.5]{suslin1991}, which is also needed
in the proof of the first exact sequence of Theorem \ref{gl-gm3}. 

\begin{lem}\label{b-h}
Let $u \in H_2(T_2)^\sigma$ and let $h \in B_2(T_2)_{T_2}$ be a representing cycle for $u$.
Let $\tau$ be the automorphism of $B_\bullet(T_2)_{T_2}$ induced by $\sigma$ and let 
$\tau(h)-h=\partial_3^{T_2}(b)$, $b \in B_3(T_2)_{T_2}$. Then the image of $u$ under the map
\[
H_2(T_2)^\sigma \arr H_3(\GM_2(R))/H_3(T_2)
\] 
coincides with the homology class of the cycle $b -\rho_s(h)$, where $s= {\mtx 0 1 1 0}$.
\end{lem}
\begin{proof}
The $E^1$-terms of the spectral sequence $E_{p,q}^2$ are of the form
\[
E_{p,q}^1=H_q(B_p(\Sigma_2) \otimes_{\Sigma_2} B_\bullet(\GM_2(R))_{T_2}).
\]
Let $h=[g_1|g_2]-[g_2|g_1]$. Then
$\tau(h)=[sg_1s^{-1}|sg_2s^{-1}]-[sg_2s^{-1}|sg_1s^{-1}]$.
By considering the commutative diagram
\[
\begin{CD}
\z \otimes_{\Sigma_2} B_3(\GM_2)_{T_2} &\leftarrow &
B_0(\Sigma_2) \otimes_{\Sigma_2} B_3(\GM_2)_{T_2} &\leftarrow&
B_1(\Sigma_2) \otimes_{\Sigma_2} B_3(\GM_2)_{T_2} \\
 @VVV  @VVV  @VVV&\\
\z \otimes_{\Sigma_2} B_2(\GM_2)_{T_2} &\leftarrow&
B_0(\Sigma_2) \otimes_{\Sigma_2} B_2(\GM_2)_{T_2} &\leftarrow&
B_1(\Sigma_2) \otimes_{\Sigma_2} B_2(\GM_2)_{T_2},
\end{CD}
\]
the necessary computations can be collected in the following diagram
\[
\begin{array}{ccccc}
b -\rho_s(h) \!\!\!\!& \leftarrow\!\!\!-\!\!\!-\!\!\!-\!\!\!-\!\!\!\shortmid &
\!\!\! s[\ ]\otimes b-[\ ]\otimes \rho_s(h) & & \\
& &
\begin{array}{c}
\vspace{-2.5mm}
-\\
\vspace{-1mm}
|\\
\downarrow
\end{array}
&& \\
&& s[\ ]\otimes h-[\ ]\otimes h & \!\!\!\!\!\!\!\!
\leftarrow\!\!\!-\!\!\!-\!\!\!-\!\!\!-\!\!\!\shortmid  & \!\!\! [s]\otimes h.
\end{array}
\]
\end{proof}

Let $R$ be a domain and $\mu_n(R)=\lan \xi \ran$, where $\xi \in R$ is a primitive 
$n$-th root of unity. Let  $\lan \xi , n, \xi \ran\in \tors(\mu(R),\mu(R))$ be the image of
$\xi$ under the following composition
\[
\mu_n(R) \overset{\simeq}{\larr} \tors(\mu_n(R),\mu_n(R))
\hookrightarrow \tors(\mu(R),\mu(R)).
\]
\begin{lem}\label{canonical}
Let $R$ be a  domain. Then we have the canonical decomposition
\[
\begin{array}{c}
H_3(T_2)= \bigoplus_{i=0}^3 H_i(\fff) \otimes_\z H_j(\fff)
\oplus \tors(\mu(R),\mu(R)),
\end{array}
\]
where a splitting map $\tors(\mu(R),\mu(R))\arr H_3(T_2)$
is defined by the formula $\lan \xi , n, \xi \ran \mapsto \chi(\xi)$, where
\begin{align*}
\chi(\xi):= \sum_{i=1}^{n}\Big(
&[(\xi,1)|(1,\xi)|(1,\xi^i)] -[(1,\xi)|(\xi,1)|(1,\xi^i)]+ \![(1,\xi)|(1,\xi^i)|(\xi,1)]\\
&\!\!\!\!\!\!\!+[(\xi,1)|(\xi^i,1)|(1,\xi)] - \![(\xi,1)|(1,\xi)|(\xi^i,1)]+[(1,\xi)|(\xi,1)|(\xi^i,1)] \Big),
\end{align*}
with $\xi \in R$ being a primitive $n$-th root of unity.
\end{lem}
\begin{proof}
For a proof see \cite[Section 4]{mirzaii-2008}.
\end{proof}

Since $R$ is a domain, $\mu(R)$ is direct limit of its finite cyclic subgroups
with a directed index set. So $\Sigma_2$ acts trivially on $\tors(\mu(R),\mu(R))$. Moreover
\begin{align*}
H_1(\Sigma_2, \mu_{2^\infty}(R)\otimes_\z \mu_{2^\infty}(R)) &
\simeq \frac{(\mu_{2^\infty}(R)\otimes_\z\mu_{2^\infty}(R))^\sigma}{(1+\sigma)
(\mu_{2^\infty}(R)\otimes_\z \mu_{2^\infty}(R))}\\
&\simeq (\mu_{2^\infty}(R)\otimes_\z\mu_{2^\infty}(R))^\sigma \\
&={}_2(\mu_{2^\infty}(R)\otimes_\z\mu_{2^\infty}(R))\\
& \simeq
\begin{cases}
0 & \text{if $\mu_{2^\infty}(R)$ is infinite or $\char(R)=\! 2$} \\
\z/2 &   \text{if $\mu_{2^\infty}(R)$ is finite and $\char(R)\neq 2$.}
\end{cases}
\end{align*}
Hence if $\mu_{2^\infty}(R)$ is infinite or $\char(R)=2$, then $\tors(\mu(R),\mu(R)) \simeq T_R$
and if $\char(R)\neq 2$ and $\mu_{2^\infty}(R)$ is finite, then we have the exact sequence
\begin{equation}\label{TR}
0 \arr \tors(\mu(R),\mu(R)) \arr T_R \arr \z/2 \arr 0.
\end{equation}
On the other hand we have
\begin{align*}
\exts(\z/2,\tors(\mu(R),\mu(R))) 
\!\simeq\!
\begin{cases} 
0 & \!\!\!\text{if $\mu_{2^\infty}(R)$ is infinite or $\char(R) = \!2$} \\
\z/2    &  \!\!\! \text{if $\mu_{2^\infty}(R)$ is finite and $\char(R)\neq 2$.}
\end{cases}
\end{align*}

\begin{lem}\label{nonsplit}
Let $R$ be a domain such that $\char(R)\neq 2$ and $\mu_{2^\infty}(R)$ is finite. 
Then the exact sequence $(\ref{TR})$ does not split, so $T_R$ is the unique non-trivial
extension of $\z/2$ by $\tors(\mu(R),\mu(R))$.
\end{lem}
\begin{proof}
The exact sequence (\ref{TR}) is, in fact, the exact sequence
\[
0 \arr H_3(T_2)/M \arr F_2H_3(\GM_2(R))/M \arr F_2H_3(\GM_2(R))/H_3(T_2) \arr 0.
\]
Let $\mu_{2^\infty}(R)=\mu_n(R)=\lan \xi \ran$, $n=2m=2^r$. Under the inclusion
\[
\{0, (-1) \otimes \xi\}={}_2(\mu_{2^\infty}(R) \otimes_\z \mu_{2^\infty}(R))
\hookrightarrow H_2(T_2)^\sigma,
\]
the image of the element $(-1)\otimes \xi$ in $H_2(T_2)^\sigma$ is represented by the cycle
$h:=[(-1,1)|(1, \xi)]-[(1,\xi)|(-1,1)]\in B_2(T_2)_{T_2}$. Now $\partial_3^{T_2}(b)=\tau(h)-h$, where
\begin{align*}
\!\!\!\!b:= & [(1,\!-1)|(1,\!-1)|(\xi, 1)]-[(1,\!-1)|(\xi,1)|(1,\!-1)]+ [(\xi,1)|(1,\!-1)|(1,\!-1)]\\
&\!\!\!\!\!\!\!\!\!\!+\!\!
\begin{array}{c}
\sum_{i=1}^{m-1}
\end{array}
\!\!\!\Big([(\xi,1)|(1,\xi)|(1,\xi^i)]-
[(1,\xi)|(\xi,1)|(1,\xi^i)]+[(1,\xi)|(1,\xi^i)|(\xi,1)]\\
& \ \ \ \ \ + [(\xi,1)|(\xi^i,1)|(1,\xi)]-[(\xi,1)|(1,\xi)|(\xi^i,1)]
+[(1,\xi)|(\xi,1)|(\xi^i,1)]\Big).
\end{align*}
Hence by Lemma \ref{b-h}, the cycle
$b-\rho_s(h)$ represents the image of $(-1) \otimes \xi$ in
$F_2H_3(\GM_2(R))/H_3(T_2)$. 
To show that our exact sequence does not splits,
it is sufficient to show that $2(\overline{b-\rho_s(h)})$ is equal to 
$\overline{\chi(\xi)}$ which by Lemma \ref{canonical}, is the image of
$\lan \xi, n, \xi \ran \in \tors(\mu(R), \mu(R))$ under the inclusion
$\tors(\mu(R), \mu(R))\simeq H_3(T_2)/M \hookrightarrow F_2H_3(\GM_2(R))/M$.
By a direct computation we have
\[
2(b-\rho_s(h))=\chi(\xi)+\partial_4(\eta(\xi)+\upsilon(\xi)),
\]
where
\begin{gather*}
\begin{array}{rl}
\vspace{2mm}
\eta(\xi)\!\!\!\!\!&:=\![s|(1,-1)|(1,-1)|(\xi,1)]- [s|(1,-1)|(\xi,1)|(1,-1)]\\
\vspace{2mm}
&\  +[(-1,1)|s|(\xi,1)|(1,-1)]-[(-1,1)|(1,\xi)|s|(1,-1)]\\
\vspace{2mm}
&\  +[(-1,1)|(1,\xi)|(-1,1)|s]-[(1,\xi)|s|(1,-1)|(1,-1)]\\
\vspace{2mm}
&\ +[s|(\xi,1)|(1,-1)|(1,-1)]-[(-1,1)|s|(1,-1)|(\xi,1)]\\
\vspace{2mm}
&\ +[(1,\xi)|(-1,1)|s|(1,-1)]-[(1,\xi)|(-1,1)|(-1,1)|s]\\
&\ +[(-1,1)|(-1,1)|s|(\xi,1)]-[(-1,1)|(-1,1)|(1,\xi)|s],
\end{array}
\end{gather*}
\begin{align*}
\upsilon(\xi):=
\!\!\!\begin{array}{c}
\sum_{i=1}^{m-1}
\end{array}
&\!\!\!\Big(
[(\xi,1)|(1,\xi)|(1,\xi^i)|(1,-1)]-[(1,\xi)|(\xi, 1)|(1,\xi^i)|(1,-1)]\\
& \!\!\!\!\!\!\!\!-[(1,\xi)|(1,\xi^i)|(1,-1)|(\xi,1)]+ [(1,\xi)|(1,\xi^i)|(\xi, 1)|(1,-1)]\\
& \!\!\!\!\!\!\!\!+[(1,\xi)|(\xi, 1)|(\xi^i, 1)|(-1,1)]-[(\xi,1)|(1,\xi)|(\xi^i,1)|(-1,1)]\\
& \!\!\!\!\!\!\!\!-[(\xi,1)|(\xi^i, 1)|(-1,1)|(1,\xi)]+[(\xi,1)|(\xi^i, 1)|(1,\xi)|(-1,1)]\Big).
\end{align*}
Thus $(\ref{TR})$ does not split.
\end{proof}

Let $n$ be a positive integer such that $2|n$. Then we have the following two non-split exact sequences
\[
\begin{CD}
0 \larr \z/n @>{\bar{r} \mapsto \overline{2r}}>> \z/2n  @>{\bar{a} \mapsto \bar{a}}>>
\z/2 \larr 0,
\end{CD}
\]
\[
\begin{CD}
0 \larr \z/2 @>{\bar{1} \mapsto \overline{n}}>> \z/2n  @>{\bar{a} \mapsto \bar{a}}>>
\z/n \larr 0.
\end{CD}
\]
Since
\[
\exts(\z/2,\z/n)\simeq
\begin{cases} 0 & \text{if $2\nmid n$} \\
\z/2    &   \text{if $2\mid n$}
\end{cases}, 
\ \ \
\exts(\z/n,\z/2)\simeq
\begin{cases} 0 & \text{if $2\nmid n$} \\
\z/2    &   \text{if $2\mid n$}
\end{cases},
\]
the first exact sequence in above is
the only nontrivial extension of $\z/2$ by $\z/n$ and the second one
is the only nontrivial extension of $\z/n$ by $\z/2$.
Let $\tors(\mu(R),\mu(R))^\sim$ denotes the unique  non-trivial extension of 
$\tors(\mu(R),\mu(R))$ by $\z/2$.

\begin{cor}\label{tr11}
For a domain $R$, $T_R\simeq \tors(\mu(R),\mu(R))^\sim$.
\end{cor}
\begin{proof}
The above observation together with Lemma \ref{TR-exact} show that when
$\char(R)\neq 2$ and $\mu_{2^\infty}(R)$ is finite, then $T_R$ can be considered 
as the unique nontrivial extension of $\tors(\mu(R),\mu(R))$ by $\z/2$. 
Therefore we have $T_R\simeq \tors(\mu(R),\mu(R))^\sim$.
Moreover if $\mu_{2^\infty}(R)$ is infinite or $\char(R)=2$, then 
$T_R\simeq \tors(\mu(R),\mu(R))\simeq \tors(\mu(R),\mu(R))^\sim$.
\end{proof}

%%%%%%%%%%%%%%%%%%%%%%%%%%%%%%%%%%%%%%%%%%%%%%%%%%%%%%%%%%%%%%%%%%%%%%%%%%
\section{The bloch-wigner exact sequence}
%%%%%%%%%%%%%%%%%%%%%%%%%%%%%%%%%%%%%%%%%%%%%%%%%%%%%%%%%%%%%%%%%%%%%%%%%%

Now we are ready to prove our main Theorem.

\begin{thm}\label{mirzaii-mokari}
Let $R$ be a ring with many units. Then we have the exact sequence
$$T_R \arr K_3^\ind(R) \arr B(R) \arr 0,$$
where $T_R$ sits in the short exact sequence
\[
0 \arr \tors(\mu(R),\mu(R))_\sigma \arr T_R \arr
H_1(\Sigma_2, \mu_{2^\infty}(R)\otimes_\z\mu_{2^\infty}(R)) \arr 0.
\]
Moreover if $R$ is an integral domain, then we have the exact sequence
\[
0 \arr \tors(\mu(R), \mu(R))^\sim \arr K_3^\ind(R) \arr B(R) \arr 0,
\]
where the composition $\tors(\mu(R),\mu(R)) \arr \tors(\mu(R), \mu(R))^\sim \arr K_3^\ind(R)$ is 
induced by the map $\mu(R) \arr \SL_2(R)$,  $\xi \mt \diag(\xi, \xi^{-1})$.
\end{thm}
\begin{proof}
Since $\im (H_3(T_3)\arr H_3(\GL_3(R)))$ is equal to 
$\im(H_3(T_2) \oplus (\rr)^{\otimes 3}\arr H_3(\GL_3(R)))$,
from Theorem \ref{gl-gm3} we get the exact sequence
\[
H_3(\GM_2(R))\oplus (\rr)^{\otimes 3} \arr H_3(\GL_3(R)) \arr B(R) \arr 0.
\]
If $B_2:= {\mtx {\rr} R 0 {\rr}}$, then for any $n\geq 0$,
$H_n(T_2)\simeq H_n(B_2)$ \cite[Theorem~1.9]{suslin1985}, \cite[Theorem 2.2.2]{guin1989}. 
The matrix $s={\mtx 0 1 1 0}$ is similar to the matrix ${\mtx 1 1 0 {-1}} \in B_2$, and hence
$\im(H_3(\Sigma_2)) \se \im(H_3(B_2))=\im(H_3(T_2))$.
Thus the above exact sequence together with (\ref{GM2}) imply the exact sequence
$F_2H_3(\GM_2(R)) \arr H_3(\GL_3(R))/N \arr B(R) \arr 0$,
where $N=\rr \cup \rr \cup \rr =\im((\rr)^{\otimes 3})$.
Now from the commutative diagram
\[
\begin{CD}
M & \hspace{-6 mm} \overset{=}{-\!\!\!-\!\!\!-\!\!\!-\!\!\!-\!\!\!\larr} & M  & \\
@VVV  @VVjV        &        \\
F_2H_3(\GM_2(R)) & \!\!\!\larr &\!\! H_3(\GL_3(R))/N & \arr  & B(R) & \arr 0,
\end{CD}
\]
we obtain the exact sequence
\begin{equation}\label{bloch1}
T_R \arr H_3(\GL_3(R))/(j(M)+N) \arr B(R) \arr 0.
\end{equation}
By an easy analysis of the Lyndon-Hochschild-Serre spectral  sequence associated to
the extension $1 \arr \SL(R) \arr \GL(R) \arr \rr \arr 1$, we obtain the exact sequence 
\[
0 \!\arr\! H_3(\SL(R))\! \arr\! H_3(\GL(R))/H_3(\GL_1(R))
\!\arr \!H_1(\rr,H_2(\SL(R))) \!\arr\! 0.  
\]
Since the action of $\rr$ on $H_i(\SL(R))$ is trivial we have 
\[
H_1(\rr,H_2(\SL(R)))\simeq \rr\otimes_\z H_2(\SL(R))\simeq \rr \otimes_\z K_2^M(R).
\]
Also this exact sequence splits, thus we have the exact sequence
\[
0 \arr \rr \otimes_\z K_2^M(R) \overset{\eta}{\arr} H_3(\GL(R))/H_3(\GL_1(R)) 
\overset{\tau}{\arr} H_3(\SL(R)) \arr 0, 
\]
where $\eta$ is given by 
$a \otimes \{b, c\} \mt {\bf c}(\diag(a,1,1), \diag(1,b,1), \diag(1, c, c^{-1}))$
and $\tau$ is induced by $A \mt \diag(\det(A)^{-1}, A)$. We have
\begin{align*}
\tau(a\cup {\bf c}(b,c)) &=\tau({\bf c}(\diag(a,1), \diag(1,b), \diag(1, c))\\
& ={\bf c}(\diag(a^{-1},a,1), \diag(b^{-1},1,b), \diag(c^{-1}, 1, c))\\
&=[b,a, c],\\
\tau(a\cup b\cup c) &= \tau({\bf c}(\diag(a,1,1), \diag(1,b,1), \diag(1, 1, c))\\
&= {\bf c}(\diag(a^{-1},a,1,1), \diag(b^{-1},1,b,1), \diag(c^{-1},1, 1, c))\\
&=[a,c,b].
\end{align*}
So we get the surjective map
$\bar{\tau}:H_3(\GL_3(R))/(j(M)+N)\two H_3(\SL(R))/T$.
Since $\ker(\tau)=\im(\eta)=\eta(\rr \otimes_\z K_2^M(R))\se j(M)+N$, we get
\begin{equation*}\label{exact2}
H_3(\GL_3(R))/(j(M)+N) \simeq H_3(\SL(R))/T.
\end{equation*}
Therefore the exact sequence (\ref{bloch1}) together with  Corollary \ref{k3-h3} and the 
exact sequence (\ref{TR-exact}) imply the first claim.

Now let $R$ be a domain. Then
from the first part and Corollary \ref{tr11} we get the exact sequence
\[
\tors(\mu(R),\mu(R))^\sim \arr K_3^\ind(R) \arr B(R) \arr 0.
\]
The natural map $\inc: \SL_2(R) \arr \GL_2(R)$ induces the commutative diagram
\[
\begin{CD}
H_3(\rr) & @>{H_3(\delta')}>> & H_3(\SL_2(R))\\
          @VVH_3(\Delta')V & & @V{H_3(\inc)}VV  \\
H_3(T_2) & @>>> & H_3(\GL_2(R)),
\end{CD}
\]
where both $\delta':\rr \arr \SL_2(R)$ and $\Delta':\rr \arr T_2$
are given by the formula $a \mt \diag(a, a^{-1})$.
Consider the commutative diagram with exact rows
\[
\begin{CD}
\!\!\!\!\!\!\!\!\!0 & \!\!\!\!\!\!\!\!\!\!\!\!\!\!-\!\!\!-\!\!\!\larr\bigwedge_\z^3\rr & 
\!\!\!\!\!\!\!\!-\!\!\!-\!\!\!-\!\!\!\larr & 
H_3(\rr) & &\!\!\!\!\!\!\!\!\!\!\!\!\!\!\!\!\!\!\!\!\!\!\!\!\!\!\!\!\!\!\!\!\!\!\!\!-\!\!\!-\!\!\!
-\!\!\!-\!\!\!-\!\!\!\larr  \tors(\rr, \rr) & \!\!\!\!\!\!\!\!\!\!\!\!\!\!\!\!\!\!\!\!\!\!\!\!\!
-\!\!\!-\!\!\!-\!\!\!-\!\!\!-\!\!\!-\!\!\!\larr 0 \\
&  @VVV  @VV{H_3(\Delta')}V @VVV & \\
0 \arr & \bigwedge_\z^3(\rr \times\rr) & \arr & H_3(\rr \times \rr)  & \ \arr\ & 
\tors(\rr \times\rr, \rr \times\rr)^{-\sigma}& \arr 0,
\end{CD}
\]
(see \cite[Lemma~5.5]{suslin1991} for the exactness of the rows). 
Let $R_1^\times=\rr \times 1$ and $R_2^\times=1 \times \rr$. Then we have
\begin{align*}
\tors(\rr \times\rr, \rr \times\rr)^{-\sigma}& =\tors(R_1^\times, R_1^\times)\oplus
\tors(R_2^\times, R_2^\times)\\
& \oplus  \Big(\tors(R_1^\times, R_2^\times)\oplus
\tors(R_2^\times, R_1^\times)\Big)^{-\sigma},
\end{align*}
which clearly from it we obtain the isomorphisms
\[
\tors(\rr, \rr)\!\overset{\simeq}{\larr}\! \Big(\tors(R_1^\times, R_2^\times)
\oplus\tors(R_2^\times, R_1^\times)\Big)^{-\sigma}\!\!\overset{\simeq}{\larr}
\!H_2(\rr \!\!\times\! \rr)/M.
\]
We denote this composition by $\varphi$. Now the commutative diagram
\[
\begin{CD}
\tors(\rr, \rr)&@>\simeq>> & H_3(\rr)/\bigwedge_\z^3\rr & 
@>>> & H_3(\SL_2(R))/j'(M')\\
@V{\simeq}VV & &@V{\overline{H_3(\Delta')}}VV & &@VVV  \\
\tors(\rr, \rr)&@>{\simeq}>\varphi> & H_3(\rr\times \rr)/M & @>>> & H_3(\GL_2(R))/j(M),
\end{CD}
\]
with $j'(M')=\im(\bigwedge_\z^3\rr \arr H_3(\SL_2(R))$, 
shows that the homomorphism $\tors(\mu(R), \mu(R)) \arr H_3(\SL(R))/T\simeq K_3^\ind(R)$
is induced by $\rr \arr \SL_2(R)$, $\xi \mt \diag(\xi, \xi^{-1})$.

To finish the proof of the theorem we have to prove that the map
\[
\tors(\mu(R), \mu(R))^\sim \arr K_3^\ind(R)
\]
is injective. Let $F$ be the quotient field of $R$ and $\bar{F}$ be
its algebraic closure. Then from the commutative diagram
\[
\begin{CD}
0 \arr & \z/2 & \arr & \tors(\mu(R),\mu(R))^\sim & \arr & \tors(\mu(R), \mu(R)) &\arr 0 \\
        &  @VV=V  @VVV @VV{\inc}V & \\
0 \arr & \z/2 & \arr & \tors(\mu(\bar{F}),\mu(\bar{F}))  & \overset{2.}{\larr}
& \tors(\mu(\bar{F}), \mu(\bar{F}))& \arr 0,
\end{CD}
\]
we see that $\tors(\mu(R),\mu(R))^\sim \arr \tors(\mu(\bar{F}),\mu(\bar{F}))$ 
is injective. Hence it is sufficient to prove the claim for $\bar{F}$. Moreover
since $K_3^M(\bar{F})$ is uniquely 
divisible \cite[Proposition 1.2]{bass-tate1973}, the torsion subgroup of $K_3(\bar{F})$ 
and $K_3^\ind(\bar{F})$ are equal.
Thus it is sufficient to prove that 
\[
\theta_{\bar{F}}: T_{\bar{F}}=\tors(\mu(\bar{F}),\mu(\bar{F})) \arr  H_3(\SL(\bar{F})\simeq 
K_3(\bar{F})
\]
is injective. From now on we assume that $F$ is algebraically close and so $F=\bar{F}$.

The injectivity of $\theta_{\C}$ is a classical result and can be proved using
Cheeger-Chern-Simon invariant $\hat{C}: H_3(\SL_2(\C)) \arr \C/\z$ \cite[\S2]{dupont1987}. 
In fact it is well-known that the composition 
\[
\q/\z\simeq \tors(\mu(\C),\mu(\C)) \overset{\theta_{\C}}{\larr} H_3(\SL_2(\C)) 
\overset{\hat{C}}{\larr} \C/\z
\] 
is injective and thus $\theta_{\C}$ is injective \cite[App. A]{dupont-sah1982}.
On the other hand, since $\mu(\bar{\q})=\mu(\C)$, by the commutative diagram
\[
\begin{CD}
 & \tors(\mu(\bar{\q}), \mu(\bar{\q})) & \overset{\theta_{\bar{\q}}}{\larr} & K_3(\bar{\q})\\
        &  @VV=V  @VVV & \\
0 \arr & \tors(\mu(\C), \mu(\C)) & \overset{\theta_{\C}}{\larr} & K_3(\C),
\end{CD}
\]
$\theta_{\bar{\q}}$ also is injective.
If $\char(F)=0$, then $\bar{\q} \se F$ and $\mu(\bar{\q})=\mu(F)$.
But Suslin has proved that for an extension $F_1 \se F_2$ of algebraically closed fields,
$K_n(F_i)$ are divisible and their torsion subgroup are isomorphic 
\cite[Main~Theorem]{suslin1983}, \cite[Proposition~3.12, Corollary~3.13]{suslin1984}. 
This implies that in the diagram 
\[
\begin{CD}
0 \arr & \tors(\mu(\bar{\q}), \mu(\bar{\q})) & \overset{\theta_{\bar{\q}}}{\larr} & 
K_3(\bar{\q})\\
&  @VV=V  @VV{\inc_\ast}V & \\
& \tors(\mu(F), \mu(F)) & \overset{\theta_{F}}{\larr} & K_3(F)
\end{CD}
\]
$\inc_\ast\circ \theta_{\bar{\q}}$ is injective. Thus $\theta_{F}$ also is injective.
Now let $\char(F)=p\neq 0$. Then $\overline{\F}_p \se F$ and thus,
with an argument similar to above, it is sufficient to prove the claim for $\overline{\F}_p$.
So we may assume that $F=\overline{\F}_p$.
Let $A$ be a henselian local domain such that its residue field 
is $F$ and its quotient field is a subfield of $\bar{\q}$ and consider the commutative diagram
\[
\begin{CD}
\tors(\mu(\bar{\q}),\mu(\bar{\q}))&@<{\psi_1}<< & \tors(\mu(A),\mu(A))&
@>{\psi_2}>>&\tors(\mu(F), \mu(F)) \\
          @V{\theta_{\bar{\q}}}VV && @V{\theta_{A}}VV & &@VV\theta_{F}V  \\
H_3(\SL(\bar{\q}))& @<{\phi_1}<< & H_3(\SL(A))& @>{\phi_2}>>&H_3(\SL(F)) \\
          @V{\simeq}VV && @V{\simeq}VV & &@VV{\simeq}V  \\
K_3(\bar{\q})& @<{}<< & K_3(A)/l(-1)K_2(A)& @>{}>>&K_3(F).
\end{CD}
\]
We have proved that $\theta_{\bar{\q}}$ is injective. Since $\psi_1$ is injective,
$\theta_{A}$ also is injective. Since $A$ is henselian, $\psi_2$ is surjective.
By Rigidity theorem of Suslin and Gabber \cite{suslin1983}, \cite{gabber1992}
if $p\nmid m$ and $n\geq 1$, then $K_n(A, \z/m)\simeq K_n(F, \z/m)$.
Since $K_n(F)$ is divisible, the commutative diagram with exact rows
\[
\begin{CD}
0 \arr & K_4(A)\otimes_\z \z/m & \arr & K_4(A, \z/m) & \arr & {}_mK_3(A) &\arr 0 \\
        &  @VVV  @VVV @VVV & \\
0 \arr & K_4(F)\otimes_\z \z/m & \arr & K_4(F, \z/m) & \arr & {}_mK_3(F)&\arr 0
\end{CD}
\]
implies that ${}_mK_3(A)\simeq {}_mK_3(F)$. Thus the map $K_3(A) \arr K_3(F)$
gives an isomorphism on torsion prime to $p$. Since $\psi_2$ is surjective,
to prove the injectivity of $\theta_{F}$, it is sufficient to prove that 
$\phi_2\circ \theta_{A}$ is injective on torsion prime to $p$. 
Let $x \in \tors(\mu(A),\mu(A))$ such that $mx=0$, $p\nmid m$, and let 
$\phi_2\circ \theta_{A}(x)=0$.
Then $\theta_A(x) \in {}_mH_3(\SL(A))$. Let $y \in K_3(A)$
such that $h_3'(y)=\theta_A(x)$. Then clearly ${2m}y=0$ and thus $2y \in {}_mK_3(A)$.
But under the isomorphism ${}_mK_3(A) \overset{\simeq}{\larr} {}_mK_3(F)$, $2y$ goes
to $2\phi_2\circ \theta_{A}(x)=0$. Thus $2y=0$, which implies that 
$2\theta(x)=h_3'(2y)=0$. If $p=2$, then from $(2,m)=1$, we see that 
$\theta_{A}(x)=0$. If $p\neq 2$, then the isomorphism
${}_2K_3(A) \overset{\simeq}{\larr} {}_2K_3(F)$ implies that $y=0$ and thus 
$\theta_A(x)=h_3'(y)=0$. In any case $\theta_A(x)=0$ and thus the injectivity of $\theta_A$
implies that $x=0$. This show that $\phi_2\circ \theta_{A}$ is injective on torsion prime to $p$
and thus the proof of the theorem is complete.
\end{proof}

\begin{rem}
(i) Our general approach to the Bloch-Wigner exact sequence
is very close to the one in \cite{mirzaii2011}.
There we introduced the group $\H_3(\SL_2(R)):=H_3(\GL_2(R))/j(M)$, which was our possible 
candidate for $K_3(R)^\ind$. Now it is easy to see that the natural map 
\[
\H_3(\SL_2(R))\arr H_3(\GL_3(R))/(j(M)+N)\simeq K_3^\ind(R)
\]
is surjective. Moreover with the techniques developed here 
it is not difficult to show that when $R$ is a domain, the
map $\H_3(\SL_2(R)) \arr K_3^\ind(R)$ is bijective (see the 
proof of Proposition \ref{SL-Kind}).
\par (ii)
The idea of the proof of the injectivity of $\tors(\mu(F), \mu(F))\arr K_3^\ind(F)$ 
for an algebraically closed field $F$, is taken from \cite[\S7]{lichenbaum1987}.
\par (iii)
We strongly believe that the map $T_R \arr K_3^\ind(R)$ in Theorem \ref{mirzaii-mokari}
is injective in general. As we have seen, this is difficult to prove even for
algebraically closed fields. In the above theorem we have shown that this is true if $R$ is 
a domain.
\end{rem}

\begin{cor}\label{hensel-local}
Let $R$ be a local ring with infinite residue field $F=R/\mmm_R$ such that the natural map
$\mu(R)\arr \mu(F)$ is injective. Then we have the exact sequence
\[
0 \arr \tors(\mu(R), \mu(R))^\sim \arr K_3^\ind(R) \arr B(R) \arr 0.
\]
\end{cor}
\begin{proof}
Since $\mu(R)\arr \mu(F)$ is injective, from
Theorem \ref{mirzaii-mokari} we obtain the exact sequence
$\tors(\mu(R), \mu(R))^\sim \arr K_3^\ind(R) \arr B(R) \arr 0$.
Now the injectivity of the left hand side map follows from the commutative diagram
\[
\begin{CD}
\ \ \ \ \ \tors(\mu(R), \mu(R))^\sim & \larr & K_3^\ind(R)\\
@VVV & \!\!\!\!\!\!\!\!\!\!\!\!\!\!\!\! @VVV \\
0 \arr \tors(\mu(F), \mu(F))^\sim & \larr & K_3^\ind(F).
\end{CD}
\]
\end{proof}

\begin{cor}\label{R-K}
Let $R$ be a discrete valuation ring with field of fraction $K$ and infinite 
residue field $F$. Suppose that $\char(K)=\char(F)$. Then the natural map
$B(R) \arr B(K)$ is an isomorphism.
\end{cor}
\begin{proof}
By \cite[Theorem~2.1]{hutchinson2014}, $K_3^\ind(R)\simeq K_3^\ind(K)$. Hence the claim follows 
from the Bloch-Wigner exact sequences for $R$ and $K$, obtained from 
Theorem \ref{mirzaii-mokari}, and the fact that $\mu(R)=\mu(K)$.
\end{proof}

\begin{exa}
(i) If $R$ is a henselian local ring with infinite residue field $F$
such that $\mu_{p^\infty}(R)=1$, where $p=\char(F)$, then $\mu(R)\arr \mu(F)$ is an isomorphism.
Thus, by Corollary \ref{hensel-local}, we have the Bloch-Wigner exact sequence over $R$.

(ii) Let $F$ be an infinite field and let $R$ be a semi-local domain containing $F$ with quotient filed
$F(X)$. Then $K_3^\ind(F)\simeq K_3^\ind(R)$ \cite[pp.~327--328]{levine1989}. 
Now as in the proof of Corollary \ref{R-K}, we can show that $B(F)\simeq B(R)$.
\end{exa}

Let $R[\epsilon]$ be the ring of dual numbers, where R is a domain. Then for any prime number 
$p$, we have
\[
\mu_{p^\infty}(R[\epsilon])\simeq
\begin{cases} 
\mu_{p^\infty}(R)    &   \text{if $p\neq \char(R)$}\\
R  &   \text{if $p= \char(R)$}
\end{cases},
\]
where $R$ is the additive group $(R,+)$. (Note that if $p=\char(R)$, then 
$\mu_{p^\infty}(R)=1$.) Since for any abelian torsion group $A$,
$A\simeq \bigoplus_{q \ prime} A_{q^\infty}$,
we have
\[
\mu(R[\epsilon])\simeq
\begin{cases} 
\mu(R)    &   \text{if $ \char(R)=0$}\\
\mu(R) \oplus R  &   \text{if $\char(R)\neq 0$.}
\end{cases}
\]
Hence if $\char(R)=0$, then $T_{R[\epsilon]}\simeq T_R\simeq \tors(\mu(R),\mu(R))^\sim$. 
If $\char(R)> 2$, then $\mu_{2^\infty}(R[\epsilon])\simeq \mu_{2^\infty}(R)$ and
\[
\tors(\mu(R[\epsilon]),\mu(R[\epsilon]))_\sigma=\tors(\mu(R), \mu(R))\oplus 
\tors(R, R)_\sigma.
\]
Moreover the map $R[\epsilon] \arr R$ given by $a+ b\epsilon \mt a$, induces a
map $T_{R[\epsilon]} \arr T_R$, which from it we get the isomorphism
\[
T_{R[\epsilon]}\simeq \tors(\mu(R),\mu(R))^\sim \oplus \tors(R, R)_\sigma.
\]
If $\char(R)=2$, then $\mu_{2^\infty}(R[\epsilon])\simeq R$ and thus
\[
H_1(\Sigma_2, \mu_{2^\infty}(R[\epsilon])\otimes_\z\mu_{2^\infty}(R[\epsilon]))\simeq 
H_1(\Sigma_2, R\otimes_\z R)\simeq R.
\]
Moreover, as in above, we have 
\[
\tors(\mu(R[\epsilon]),\mu(R[\epsilon]))_\sigma=\tors(\mu(R), \mu(R))\oplus 
\tors(R, R)_\sigma.
\]
These imply that $T_{R[\epsilon]} \simeq \tors(\mu(R),\mu(R))^\sim \oplus T_\epsilon$,
where  $T_\epsilon$ sits in an exact sequence
$0 \arr \tors(R, R)_\sigma \arr T_\epsilon \arr R \arr 0$.
Thus we have proved the following proposition.

\begin{prp}\label{dual-numbers}
Let $R[\epsilon]$ be the ring of dual numbers, where R is a domain with many units.
Then we have the exact sequence 
\[
\tors(\mu(R), \mu(R))^\sim \oplus T_\epsilon\arr K_3^\ind(R[\epsilon]) \arr B(R[\epsilon]) \arr 0,
\]
where $T_\epsilon$ is $0$ if $\char(R)=0$, is isomorphic to $\tors(R, R)_\sigma$ if $\char(R)> 2$
and sits in an exact sequence $0 \arr \!\tors(R, R)_\sigma\! \arr T_\epsilon\! \arr\! R \arr 0$
if $\char(R)=2$. Moreover the map $\tors(\mu(R), \mu(R))^\sim\arr K_3^\ind(R[\epsilon])$
always is injective.
\end{prp}
\begin{proof}
The ring $R[\epsilon]$ is a ring with many units \cite[Proposition 2.5]{elbaz1999}.
Now the claim follows from Theorem \ref{mirzaii-mokari} and the above computations.
\end{proof}

%%%%%%%%%%%%%%%%%%%%%%%%%%%%%%%%%%%%%%%%%%%%%%%%%%%%%%%%%%%%%%%%%%%%%%%%%%%
\section{Third homology group of $\SL_2$}
%%%%%%%%%%%%%%%%%%%%%%%%%%%%%%%%%%%%%%%%%%%%%%%%%%%%%%%%%%%%%%%%%%%%%%%%%%%

Let $\alpha: H_3(\SL_2(R))_\fff \arr K_3(R)^\ind$
be the composition of the following sequence of maps
\[
H_3(\SL_2(R))_\fff \arr H_3(\SL(R))
\overset{({\bar{h}}_3')^{-1}}{{{-\!\!\!-\!\!\!-\!\!\!\larr}}} K_3(R)/l(-1)K_2^M(R)
\twoheadrightarrow K_3(R)^\ind.
\]
It was known for very long time that $\alpha$ is an isomorphism over algebraically 
closed fields \cite[Theorem~4.1]{sah1989}. Following this, Suslin raised the following
question:
\medskip
~\\
{\bf Question} (Suslin).
Let $F$ be an infinite field. Is it true that the group $H_0(F^\times, H_3(\SL_2(F)))$ 
coincides with $K_3(F)^\ind$?
(See \cite[Question 4.4]{sah1989}).
\medskip

Hutchinson and Tao proved that $\alpha$ always is surjective 
\cite[Lemma 5.1]{hutchinson-tao2009}. They also gave some criterion for the injectivity 
of $\alpha$ (see \cite[p.~1680]{hutchinson-tao2009}). Their criterion have been improved in 
\cite[Theorem 4.4]{mirzaii2014}. We collect these result in the following proposition.

\begin{prp}\label{SL-Kind}
Let $F$ be an infinite field. Then
\par {\rm (i)} $\alpha$ always is surjective,
\par {\rm (ii)} $\alpha$ is injective if and only if the maps
$H_3(\SL_2(F))_{F^\ast} \arr H_3(\GL_2(F))$ and 
$H_3(\GL_2(F)) \arr H_3(\GL_3(F))$ are injective.
\end{prp}
\begin{proof}
For the proof of (i) see \cite[Lemma 5.1]{hutchinson-tao2009}. To prove (ii) first note that
using the exact sequence $H_3(\GM_2(F)) \arr H_3(\GL_2(F)) \arr B(F) \arr 0$
from Theorem \ref{gl-gm3} and replacing $\H_3(\SL_2(F)):=H_3(\GL_2(F))/j(M)$ with 
$H_3(\GL_3(F))/(j(M)+N)$ in the proof of Theorems \ref{mirzaii-mokari}, one obtains
the exact sequence
\[
0 \arr \tors(\mu(F), \mu(F))^\sim \arr \H_3(\SL_2(F)) \arr B(F) \arr 0.
\]
Now the rest of the proof can be proceeded as in \cite[Theorem 4.4]{mirzaii2014}.
\end{proof}

\begin{exa}
It is known that the natural maps $H_3(\GL_2(\R))  \arr  H_3(\GL_3(\R))$ and
$H_3(\SL_2(\R))_{\R^\ast}  \arr  H_3(\SL(\R))$ are injective
\cite[Theorem 5.4, Proposition 6.1]{mirzaii2008}. Since 
$H_3(\SL(\R)) \arr H_3(\GL(\R))\simeq H_3(\GL_3(\R))$ is injective,
one easily can show that $H_3(\SL_2(\R))_{\R^\ast}  \arr  H_3(\GL_2(\R))$ also is injective.
On the other hand Parry and Sah  have proved that the action of $\R^\ast$ on $H_3(\SL_2(\R))$
is trivial \cite[App.~C]{parry-sah1983}, so $H_3(\SL_2(\R))_{\R^\ast}=H_3(\SL_2(\R))$.
Now by Proposition \ref{SL-Kind}, we have
$H_3(\SL_2(\R))\simeq K_3^\ind(\R)$.
Therefore by Theorem \ref{mirzaii-mokari}, we have the Bloch-Wigner exact sequence
\[
0 \arr \z/4 \arr H_3(\SL_2(\R)) \arr B(\R) \arr 0.
\]
\end{exa}
\medskip

It is well-known that when $R$ is a ring with many units, the kernel of the maps 
$H_3(\SL_2(R))_{R^\ast} \arr H_3(\GL_2(R))$ and 
$H_3(\GL_2(R)) \arr H_3(\GL_3(R))$ consist of 2-torsion elements 
\cite[Corollary~3.3, Lemma~3.6]{mirzaii2012}.
Using these facts one can show that 
\[
H_3(\SL_2(R), \z[1/2])_{R^\ast} \simeq K_3^\ind(R)_{\z[1/2]},
\] 
\cite[Theorem 3.7]{mirzaii2012}. Moreover if $\rr=\rr^2$
we have $H_3(\SL_2(R)) \simeq K_3^\ind(R)$. Combining these facts with the results
of the previous section we obtain the following proposition, which first has been proved in 
\cite[Corollary 5.4]{mirzaii2011}. 

\begin{prp}
Let $R$ be a ring with many units. Then we have the exact sequence
\[
{\tors(\mu(R),\mu(R))_\sigma}_{\z[1/2]} \arr H_3(\SL_2(R), \z[1/2])_{R^\ast} 
\arr B(R)_{\z[1/2]} \arr 0.
\] 
If $\rr=\rr^2$, then we have the exact sequence
\[
\tors(\mu(R),\mu(R))_\sigma \arr H_3(\SL_2(R)) \arr B(R) \arr 0.
\]
Moreover if $R$ is a domain, then the left hand side maps in both exact sequences are 
injective.
\end{prp}

%%%%%%%%%%%%%%%%%%%%%%%%%%%%%%%%%%%%%%%%%%%%%%%%%%%%%%%%%%%%%%%%%%%%%%%%%%%

\bigskip
\address{{\footnotesize
Behrooz Mirzaii,
%Department of Mathematics (ICMC),
University of Sao Paulo (USP),
Sao Carlos, Brazil.

e-mail:\ bmirzaii@icmc.usp.br,
}}

\bigskip
\address{{\footnotesize
Fatemeh Y. Mokari,
%Department of Mathematics (IMECC),
University of Campinas (Unicamp),
Campinas, Brazil.

email:\ f.mokari61@gmail.com
}}
\end{document}